\documentclass[12pt,twoside]{article}
\usepackage[margin=3cm]{geometry}
\usepackage{mathtools}
\mathtoolsset{showonlyrefs}
\usepackage{float}
\usepackage{amsmath,amssymb,amsthm}
\usepackage[american]{babel}
\usepackage[latin1]{inputenc}
\usepackage[numbers]{natbib}
\usepackage{mathrsfs}
\usepackage{booktabs}
\usepackage{url}
\usepackage{amscd}
\usepackage{amstext}
\usepackage{amsfonts}
\usepackage{bbm}
\usepackage{soul}
\usepackage[latin1]{inputenc}
\usepackage[T1]{fontenc}
\usepackage{enumerate}
\numberwithin{equation}{section}

\newcommand{\mathset}[1]{\mathbbm{#1}}
\newcommand{\morf}[4][\to]{ #2 \colon #3 #1 #4}
\newcommand{\map}[4][\mapsto]{ #2 \colon #3 #1 #4}
\newcommand{\F}{\mathcal{F}}
\newcommand{\abs}[2][]{#1\lvert #2#1\rvert}
\newcommand{\norm}[2][]{#1\lVert #2#1\rVert}
\newcommand{\inner}[2][]{#1\langle #2 #1\rangle}
\newcommand{\R}{\mathset{R}}
\newcommand{\V}{\mathcal{V}}
\newcommand{\B}{\mathscr B}
\newcommand{\N}{\mathset{N}}
\newcommand{\pr}[2][]{#1 ( #2_t #1)_{t\in\R}}
\newcommand{\D}{\mathset{D}}

\newcommand{\xT}{(X_t)_{t\in T}}
\newcommand{\E}{\mathbb{E}}
\renewcommand{\P}{\mathbb{P}}
\newcommand{\1}{\mathbf{1}}

\newcommand{\dist}{\stackrel{\scriptscriptstyle\smash{\mathrm{d}}}{=}}

\def \lt{[\hskip-1.3pt[ }
\def \rt{ ]\hskip-1.3pt]}

\theoremstyle{plain}
\newtheorem{lemma}{Lemma}[section]
\newtheorem{proposition}[lemma]{Proposition}
\newtheorem{theorem}[lemma]{Theorem}
\newtheorem{corollary}[lemma]{Corollary}
\theoremstyle{definition}
\newtheorem{example}[lemma]{Example}

\newtheorem{remark}[lemma]{Remark}
\theoremstyle{definition}

\newcommand{\devnull}[1]{}

\title{Characterization of the finite  variation property for a class of stationary increment 
infinitely divisible processes}
\author{Andreas Basse-O'Connor$^{*\dagger}$ and Jan Rosi\'nski$^{*\ddagger}$
\\  
  $^*$Department of Mathematics, The University of Tennessee, USA.
  \\  $^\dagger$E-mail: basse@imf.au.dk \qquad $^\ddagger$E-mail: rosinski@math.utk.edu }
  
  \date{December 9, 2012}
  
\begin{document}
 
 \maketitle
\begin{abstract}
We characterize the finite variation  property for stationary increment mixed moving averages driven by  infinitely divisible random measures. Such processes include  fractional and moving average processes driven by L\'evy processes, and also their mixtures. We  establish two types of zero-one laws for the finite variation property. We also consider some examples  to illustrate our results. 
 
\medskip 
\noindent
\textit{Keywords: finite variation; infinitely divisible processes; stationary processes; fractional processes; zero-one laws} 

\smallskip
\noindent
\textit{AMS Subject Classification: 60G48; 60H05; 60G51; 60G17} 
\end{abstract}

%%%%%%%%%%%%%%%%%%%%%%%%%%%%%%
\section{Introduction}\label{sec_intro}
%%%%%%%%%%%%%%%%%%%%

Processes with stationary, but not necessarily independent, increments have always been  of interest in probability and its applications. They are used to model long memory phenomena. Examples include fractional and moving average processes driven by bilateral L\'evy processes, as well as their superpositions called mixed fractional and mixed moving average processes, respectively. It has been of interest to determine when such processes are semimartingales and, in particular, when they have locally finite variation. Such questions for Gaussian moving averages were resolved  by \citet[Theorem~6.5]{Knight}.  Recently, \citet{Basse_Pedersen} characterized the semimartingale and finite variation properties for stochastic convolutions of non-Gaussian L\'evy processes but their arguments do not apply to moving averages. \citet{Be_Li_Sc} gave necessary and sufficient conditions for square integrable fractional L\'evy processes to have sample paths of finite variation and showed that the total variation property, for these processes, satisfies a zero-one law. 

In this paper we characterize the finite variation property for a wide class of stationary increment infinitely divisible processes that includes fractional L\'evy processes, moving averages and mixtures of these processes. We also establish two types of zero-one laws for such processes. Therefore, we extend results of \cite{Knight} and  \cite{Be_Li_Sc} to a much larger class of processes but our methods are different.  Our work utilizes Banach space techniques, the crucial observation that $BV[0,1]$, the space of functions of finite variation, is a Banach space of cotype~2, and  arguments in the spirit of  \citet[Theorem~24]{Hardy_Littlewood}.

The paper is organized as follows. 
In Section~\ref{sec_MA} we define the class of processes we consider. They are Stationary Increment Mixed Moving Average type (SIMMA for short) processes, see \eqref{eq:mma} and \eqref{eq:ma}. In Section~\ref{main_results} we state the main results of this paper. Theorem~\ref{th:3} gives sufficient conditions for a SIMMA process to have finite variation.  Theorem~\ref{necessary_fv}, which is the most difficult result of this work, gives necessary conditions. Theorems~\ref{th:1} and \ref{th:2} state the zero-one laws. In Section~\ref{sub_ex} we determine the finite variation property on examples of processes driven by mixtures of stable random measures and tempered stable random measures. Sections~\ref{proof_fv} and \ref{proof_01} contain proofs of the main results.
\bigskip	
   
 %%%%%%%%%%%%%%%%%%%%%%%%%%%%%%%%%%%%%
\section{Preliminaries} \label{sec_MA}
%%%%%%%%%%%%%%%%%%%%%%%%%%%%%%%%%%%%%

Throughout this paper $(\Omega,\F,\P)$ stands for a probability space and
 $(V, \mathcal V, m)$ denotes a $\sigma$-finite measure space. Let $\lambda$ be the Lebesgue measure on $\R$, $\B_0 = \{B \in \B(\R): \lambda(A) < \infty\}$, and let ${\mathcal V}_0 = \{B \in \mathcal{V}: m(B)< \infty\}$.  
 Consider a stationary increment mixed moving average (SIMMA, for short) process $X=\pr X$ given by 
 \begin{equation} \label{eq:mma}
X_t=\int_{\R\times V} \big(f(t-s, v)- f_0(-s, v)\big)\, W(ds, dv),\qquad t\in\R,
\end{equation} 
where $f, f_0 : \R\times V \mapsto \R$ are measurable deterministic functions and  $W$ is an infinitely divisible independently scattered random measure (random measure, for short) defined on the $\sigma$-ring generated by $\B_0\times {\mathcal V}_0$ such that for all $A\in \B_0$,  $B \in {\mathcal V}_0$ and $u\in\R$,
\begin{align}\label{eq:M}
&\mathbb{E} e^{iu W(A\times B)} \\ 
&\quad = \exp\Big[\lambda(A) \int_{B} \Big( iu\theta(v)-\frac{1}{2}u^2 \sigma^2(v) +\int_{\R} (e^{iux}-1-iu\lt x\rt) \, \rho_v(dx)\Big)\, \,m(dv)\Big]. 
\end{align}
Here  $\rho=\{\rho_v:v\in V\}$ is a measurable family of  L\'evy measures on $\R$, so that, for each $v\in V$, $\rho_v$ is a L\'evy measure on $\R$ and for all $A\in \V$, $v\mapsto \rho_v(A)$ is measurable. The functions 
$\morf{\theta}{V}{\R}$ and  $\morf{\sigma^2}{V}{\R_+}$ are measurable,  $x\mapsto \lt x\rt=x/( |x| \vee 1)$  is a truncation function on $\R$. 

The integral in \eqref{eq:mma} is defined as in \citet[page~460]{Rosinski_spec}.   According to \cite[Theorem 2.7]{Rosinski_spec},  given a measurable deterministic function $\morf{\phi}{\R\times V}{\R}$, the integral $\int_{\R\times V} \phi(s,v) \, \Lambda(ds,dv)$ exists	if and only if
\begin{enumerate}[(a)]
   \item \label{i1} $\int_{\R\times V} |B(\phi(s,v), v)| \, ds\,m(dv) < \infty$,
      \item \label{i2} $\int_{\R\times V} K(\phi(s,v), v) \, ds\,m(dv) < \infty$,
    \end{enumerate} 
where 
\begin{align} \label{B}
B(x, v) = {}& x b(v) + \int_{\R} \big( \lt xy\rt - x \lt y\rt\big) \, \rho_{v}(dy)\quad \text{and}
\\
\label{K}
K(x, v) = {}& x^2 \sigma^2(v) + \int_{\R}  \lt xy\rt^{2} \, \rho_{v}(dy), \qquad x \in \R, \ v \in V.
\end{align}
We further assume that $W$ is purely stochastic, that is
\begin{equation}\label{eq:pur-sto}
m(v:\rho_v(\R)=0,\, \sigma^2(v)=0)=0.
\end{equation}

The process $X$ in \eqref{eq:mma} is infinitely divisible, i.e., all its finite dimensional distributions are infinitely divisible. Since $W$ is invariant in distribution under the shift on $\R$, $\pr X$ has stationary increments and thus is continuous in probability, cf.\  \cite[Proposition~2.1]{Rosinski_preprint}. When $f_0 \equiv 0$ in \eqref{eq:mma}, then $\pr X$ is a mixed moving average process (cf.\  \cite{mixed_ma_R}).  If $V$ is a one-point space, then the $v$-component can be removed from \eqref{eq:mma}--\eqref{eq:M} and $W$ becomes a random measure generated by increments of a two-sided L\'evy process that we also denote by $W$. In this case $\pr X$ is called a stationary increment moving average (SIMA) process written as
\begin{equation}\label{eq:ma}
X_t = \int_{\R} \big(f(t-s)- f_0(-s)\big)\, dW_s, \qquad t\in\R.
\end{equation}
If also $f(s)=f_0(s)=s_+^\alpha$ for some $\alpha\in \R$ and $s_+=\max\{0,s\}$, then $\pr X$ is a  linear fractional L\'evy process. If $f_0 \equiv 0$, then $\pr X$ is a moving average.  Overall, SIMMA processes cover a large class of stationary increment infinitely divisible processes of interest.

We will often consider a symmetrization $\bar{X}=\pr {\bar{X}}$ of a process $X=\pr X$ defined as $\bar{X}_t = X_t-X'_t$, where the process $X'$ is an independent copy of $X$.  If  $X$ is a SIMMA process given by  \eqref{eq:mma}, then so is  its symmeterization $\bar{X}$. In this case, $\bar{X}$ is  given by \eqref{eq:mma} with $W$ replaced by  $\bar{W}$, where $\bar{W}$ is a symmetrization of $W$ defined analogously. 

Let $I\subseteq \R$ be an interval. A function $\morf{h}{I}{\R}$ is said to be of finite variation, if for all $a,b\in I$ with $a<b$,
\begin{equation}\label{}
 \norm{h}_{BV[a,b]}:= \sup_{\stackrel{a=t_0<\cdots< t_n=b}{n \in \N}} \ \sum_{k=1}^{n} |h(t_k) - h(t_{k-1})|<\infty.
\end{equation}
For example, if $h$ is   absolutely continuous, that is,  there exists a locally integrable function $\dot h$ such that 
 \begin{equation}\label{}
 h(t)-h(u)=\int_u^t \dot h(s) \,ds, \qquad u,t\in I,\ u<t, 
\end{equation}
 then $h$ is of finite variation and $ \norm{h}_{BV[a,b]} = \int_a^b |\dot h(s)|\,ds$.
  
We will always choose a separable process $X=\pr X$ satisfying \eqref{eq:mma}. Since $X$ is continuous in probability, we may and do assume that the set $\D \subset \R$ of dyadic numbers is its separant, see \cite{Gikhman}. Then
\begin{equation}\label{eq:fvX}
\norm{X}_{BV[a,b]}=\sup_{n\in \N} \, \sum_{i=1}^{k_n} |X_{t_i^n} -X_{t_{i-1}^n}|\qquad \text{a.s.},
\end{equation}
where $a=t_0^n<\cdots< t_{k_n}^n=b$ are such that $\{t_i^n\}_{i=1}^{k_n-1} \subset \D$ and $\max_{1\leq i\leq k_n} t_i^n-t_{i-1}^n\to 0$ as $n\to \infty$. Similarly, we may view  $f_t:=f(t-\cdot,\cdot)$, $t\in \R$, as a stochastic process with respect to some probability measure $Q$ on $\R\times V$ that is equivalent to $\lambda\otimes m$. Since $X$ is continuous in probability, so is its symmetrization $\bar{X}$. It follows from  \cite[Theorem 3.4 and Proposition 3.6(i)]{Rosinski_spec}, applied to $\bar{X}$,  that  the map $t \mapsto f_t$ is continuous in $Q$-measure. Thus we may and do assume that $(f_t)_{t\in\R}$ is separable relative to probability $Q$, with $\D$ being a separant. Consequently, we have
\begin{equation}\label{eq:fvf}
\norm{f}_{BV[a,b]}=\sup_{n\in \N} \, \sum_{i=1}^{k_n} |f_{t_i^n} -f_{t_{i-1}^n}| \qquad \lambda\otimes m\text{-a.e.},
\end{equation}
where $t_i^n$ are as in \eqref{eq:fvX}.

%%%%%%%%%%%%%%%%%%%%%%
\section{Main results}\label{main_results}
%%%%%%%%%%%%%%%%%%%%%

%%%%%%%%%%%%%%%%%
\subsection{Characterization of finite variation}\label{cha_fv}
%%%%%%%%%%%%%%%

Here we give closely related sufficient and necessary conditions for SIMMA processes to have paths of finite variation. 

\begin{theorem}[Sufficiency]\label{th:3}
Let $X=\pr X$ be a process given by \eqref{eq:mma}. Suppose that for $m$-a.e.\ $v$, $f(\cdot, v)$ is absolutely continuous and its derivative $\dot f(s, v)= \frac{\partial }{\partial s}f(s,v)$ satisfies the following two conditions
  \begin{equation} \label{gen_eq}
C_f:= \int_\R \int_V \abs{\dot f(s,v)}^2 \,  \sigma^2(v)\,m(dv)\,ds<\infty,
\end{equation}
 and
 \begin{equation}\label{eq:Cf}
 D_f:=\int_\R\int_V \int_{\R} \big( |x{\dot f}(s,v)|^2 \wedge |x{\dot f}(s,v)|\big) \, \rho_v(dx) \, m(dv)\, ds < \infty.
 \end{equation}
  Then $\pr X$  has absolutely continuous sample paths a.s.\ whose total variation is integrable on each bounded interval. Moreover,  $\lambda\otimes \P$-a.e.\
\begin{equation} \label{}
\frac{dX_t}{dt} =  \int_{\R\times V} {\dot f}(t-s, v) \, W(ds, dv),\qquad t\in\R,
\end{equation}
where  the right hand side is a well-defined mixed moving average process with paths in $L^1$ a.s.\ on each finite interval.
  \end{theorem}
  
  \medskip	
    
 \begin{corollary}\label{cor_var}
  In addition to the assumptions of Theorem \ref{th:3}, suppose that $W$ is a mean zero random measure. Then
\begin{equation} \label{}
\E \|X\|_{BV[0,1]} \le  (2/\pi)^{1/2}C_f^{1/2}+(5/4)\max\{D_f, D_f^{1/2}\}. 
\end{equation}
\end{corollary}
      
\bigskip	

The converse to Theorem \ref{th:3} is more complex due to the vast class of possible random measures $W$. Assumption \eqref{inf-var}  precludes $W$ having locally finite variation, which necessitates $f$ to have absolutely continuous sections (see Remark \ref{remark_1}). 
   
  \begin{theorem}[Necessity]\label{necessary_fv}
Suppose that $X$ has paths of finite variation a.s.\ on $[0,1]$ and 
for $m$-almost every $v \in V$ we have either
\begin{equation}\label{inf-var}
  \int_{-1}^1 |x| \, \rho_{v}(dx)=\infty\quad \text{or}\quad \sigma^2(v)>0.
\end{equation}
Then for $m$-a.e.\ $v$, $f(\cdot, v)$ is absolutely continuous, its derivative $\dot f(\cdot, v)$ satisfies \eqref{gen_eq} and 
\begin{equation}\label{eq:1suf}
 \int_\R\int_\R \big(\abs{\dot f(s,v)x}\wedge \abs{\dot f(s,v)x}^2\big) (1 \wedge x^{-2}) \,\rho_v(dx)\,ds<\infty \qquad m\text{-a.e.}
  \end{equation}
If, additionally,
\begin{equation}\label{eq:u0}
\limsup_{u \to \infty} \, \frac{u\int_{\abs{x}>u} \abs{x}\,\rho_v(dx)}{\int_{|x|\le u} x^2 \, \rho_v(dx)} < \infty \qquad m\text{-a.e.}
\end{equation}
then  $\dot f(\cdot, v)$ satisfies \eqref{gen_eq} and 
\begin{equation}\label{fdot_int}
 \int_{\R} \int_{\R} ( |x{\dot f}(s,v)|^2 \wedge |x{\dot f}(s,v)|) \, \rho_v(dx)\, ds < \infty \quad m\text{-a.e.}
 \end{equation}
Finally, if 
 \begin{equation}\label{eq:u00}
\sup_{v\in V}\sup_{u > 0} \, \frac{u\int_{\abs{x}>u} \abs{x}\,\rho_v(dx)}{\int_{|x|\le u} x^2 \, \rho_v(dx)} <\infty
\end{equation}
then $\dot f(\cdot, v)$ satisfies \eqref{gen_eq} and \eqref{eq:Cf}.
\end{theorem}

Notice that \eqref{eq:u0} and \eqref{eq:u00} are well-defined by the convention $a/0:=\infty$ if $a\in [0,\infty]$.  

\begin{remark}\label{}
Theorem \ref{necessary_fv} constitutes a complete converse to Theorem \ref{th:3} when \eqref{inf-var} holds and either  \eqref{eq:u0} holds and $V$ is finite or \eqref{eq:u00} holds. 
\end{remark}

 Surprisingly, it is not easy to find a centered random measure $W$ failing  \eqref{eq:u0}. Below we will give conditions under which \eqref{eq:u0} or \eqref{eq:u00} hold.
 Recall that a measurable function $\morf{h}{\R_+}{(0,\infty)}$ is  regularly varying at $\infty$ (resp.\ at $0$) of index $\beta\in \R$ if for all $a>0$, $h(a t)/h(t)\to a^\beta$ as $t\to \infty$ (resp.\ as $t\to 0$), see \cite[page~18]{Bingham}.
A measure $\mu$ on $\R$ is said to be regularly varying if $x\mapsto \mu([-x,x]^c)$  is a regularly varying function.  

\begin{proposition}\label{remark}  
Suppose that $\rho_v(\R)>0$ for all $v\in V$. Then \eqref{eq:u0} is satisfied when one of the following two conditions holds for $m$-almost every $v \in V$
\begin{itemize}
	\item[(i)] $\int_{\abs{x}>1} x^2 \, \rho_v(dx)<\infty$ or
	\item[(ii)]  $\rho_v$ is regularly varying at $\infty$ of index $\beta\in [-2,-1)$.
\end{itemize}
Condition \eqref{eq:u00} holds when $\rho_v=\rho_0$ for all $v$ and some fixed L\'evy measure $\rho_0$ satisfying \eqref{eq:u0} and such that $\rho_0$ is regularly varying of index $\beta\in (-2,-1)$ at  0. 
 \end{proposition}

 \begin{proof}[Proof of Proposition~\ref{remark}]
   \emph{(i)}:   For $v\in V$  choose  $u_0 = u_0(v) >0$ such that   $\rho_v([-u_0,u_0])>0$. For all $u>u_0$ 
 \begin{align}\label{}
\frac{u\int_{\abs{x}>u}\abs{x}\,\rho_v(dx) }{\int_{\abs{x}\leq u} x^2\,\rho_v(dx)} \leq \frac{\int_{\abs{x}>u} x^2 \,\rho_v(dx) }{\int_{\abs{x}\leq u} x^2\,\rho_v(dx)} \leq \frac{\int_{\abs{x}>u_0} x^2 \,\rho_v(dx) }{\int_{\abs{x}\leq u_0} x^2\,\rho_v(dx)} \,,
\end{align}
which proves \emph{(i)}. 
 
 \emph{(ii)}: Set  $g(r)=\rho_v([-r,r]^c)$ for $r>0$.  For all $u>0$ we have
 \begin{align}\label{eq-ii-1}
{}&   \int_{\abs{x}>u} \abs{x}\,  \rho_v(dx)=  
 \int_0^\infty \rho_v(x: |x|>r,|x|>u)\,dr\\ {}& \qquad \label{eq-ii-2}=
  u \rho_v([-u,u]^c) +\int_u^\infty \rho_v([-r,r]^c)\,dr= u g(u)+\int_u^\infty g(r)\,dr.
 \end{align} 
Moreover, 
 \begin{align}\label{eq-ii-3} 
  {}& \int_{\abs{x}\leq u} x^2\,\rho_v(dx)= \int_0^\infty 2r\rho_v(x: |x|>r, |x|\leq u)\,dr\\
    {}&  \qquad = \int_0^u 2r\big( \rho_v([-r,r]^c)-\rho_v([-u,u]^c)\big)\,dr
  = \int_0^u 2r \rho_v([-r,r]^c)\,dr - u^2\rho_v([-u,u])\\ \label{eq-ii-4} 
  {}&\qquad  =\int_0^u 2rg(r)\,dr - u^2g(u).
 \end{align}
 Since $g$ is regularly varying at $\infty$ we may choose $u_0=u_0(v)$ such that $g(u)>0$ for all $u\geq u_0$.
 Since $g$ is locally bounded and regularly varying at $\infty$ of index $\beta$ we have by  Karamata's Theorem \cite[Theorem~1.5.11]{Bingham}  
 as $u\to\infty$ ($u>u_0$) \noeqref{eq-ii-2}\noeqref{eq-ii-3}\noeqref{eq-ii-4}
 \begin{align}\label{two-limit}
\frac{ug(u)}{\int_u^\infty g(r)\,dr}\to -(\beta+1) 
\qquad \text{and}\qquad 
\frac{u^2g(u)}{\int_ 0^u r g(r)\,dr}\to \beta+2.
 \end{align} 
 For  the first limit we have used that $\beta<-1$ and for the second limit that  $\beta\geq -2$. 
By   \eqref{eq-ii-1}--\eqref{two-limit} and  by dividing both the numerator and denominator by  $u^2 g(u)$  we have for all $u>u_0$
 \begin{equation}
 \frac{ u \int_{\abs{x}>u} \abs{x} \,\rho_v(dx)}{\int_{\abs{x}\leq u} x^2\,\rho_v(dx)} = \frac{1+ \int_u^\infty g(r)\,dr/(ug(u))}{\int_ 0^u 2r g(r)\,dr/(u^2g(u))-1}\xrightarrow[u\to \infty]{} \frac{1-(\beta+1)^{-1}}{2(\beta+2)^{-1}-1}.
 \end{equation}
 (The limit should be understood as 0 when  $\beta=-2$.) This shows 
\eqref{eq:u0}. 
 
The proof of the last part of this proposition is similar to the proof of \emph{(ii)} and is thus  omitted.  
\end{proof} 
 
  \begin{remark}\label{remark_1}
As we mentioned earlier, Condition~\eqref{inf-var} is in general necessary to deduce that $f$ has absolutely continuous sections. Indeed, let $V$ be a one point space so that $W$ is generated by increments of a L\'evy process denoted again by $W$. If \eqref{inf-var} is not satisfied, then taking  $f=\1_{[0,1]}$ we get that $X_t=W_t-W_{t-1}$ is of finite variation, but $f$ is not continuous.  
\end{remark}
  
\bigskip	

 %%%%%%%%%%%%%%%%%%
\subsection{Zero-one laws}\label{sub_01}
%%%%%%%%%%%%%%%%%%

We distinguish two types of zero-one laws, a global one which always holds and a local one holding only in certain situations.

\begin{theorem}[Global 0-1]\label{th:1}  
Let $X=\pr X$ be a process given by \eqref{eq:mma}. Then
\begin{equation} \label{eq:0-1}
\P\big(\|X\|_{BV[a,b]}< \infty \quad \text{for all } a < b \big)=0 \ \text{or} \ 1.
\end{equation}
\end{theorem}

\bigskip	

\begin{theorem}[Local 0-1]\label{th:2}
Let $a<b$ be fixed reals. Then,
\begin{equation} \label{eq:0-1a}
\P\big(\|X\|_{BV[a,b]}< \infty\big)=0 \ \text{or} \ 1
\end{equation}
provided one of the following conditions is satisfied:
\begin{enumerate}[(a)]
\item \label{local_1}
$f(\cdot,v)$ is of finite variation for  $m$-a.e.\ $v$, 
 \item\label{local_2}
$\rho_v(\R)=\infty$ for $m$-a.e.\ $v$.
\end{enumerate}
Furthermore, if $\P\big(\|X\|_{BV[a,b]}< \infty\big)=1$, then  (a) holds. 
\end{theorem}

\begin{remark}\label{exs_count}
{\rm 
The following example shows that the local zero-one law does not always hold. 
Let $[a,b]=[0,1]$ and let $\morf{f}{\R}{\R}$ be a continuous function such that $f$ has infinite total variation on each subinterval of $[0,1]$ and $f(x)=0$ for $x \in [0,1]^c$. Consider the case where there is no $v$-component and $W$ is given by \eqref{eq:M} with $\sigma^2=\theta=0$ and $\rho=\delta_1$. Then $\{W(A): A\in \B_0 \}$ is a Poisson random measure with Lebesgue intensity measure.    Consider a moving average process
\begin{equation}\label{}
  X_t =\int_{\R} f(t-s)\, W(ds).
\end{equation}
Given that $W([-1,1])=0$, \, $X_t=0$ for all $t\in [0,1]$. Thus
$$
\P( \|X\|_{BV[0,1]} < \infty) \ge \P\big(W([-1,1])=0\big) = e^{-2}.
$$
Also, given that $W([-1,0))=0$ and $W([0,1])=1$, \, there exists an $s \in [0,1]$ such that  $X_t=f(t-s)$ for all $t\in [0,1]$, which implies  that $\|X\|_{BV[0,1]}=\infty$.  Hence
\begin{align*}
& \P( \|X\|_{BV[0,1]}  < \infty)  = 1- \P( \|X\|_{BV[0,1]} =\infty) \\ & \qquad 
 \le 1-  \P\big(W([-1,0))=0, \ W([0,1])=1\big)\\ & \qquad 
 = 1-  \P\big(W([-1,0))=0\big) \P\big(W([0,1])=1\big) = 1- e^{-2}.
\end{align*}
This shows   $\P( \|X\|_{BV[0,1]} < \infty)\in [e^{-2},1-e^{-2}]$. 
\qed
}
\end{remark}

\bigskip	

%%%%%%%%%%%%%%%%%
\section{Examples}\label{sub_ex}
%%%%%%%%%%%%%%%%%

In this subsection we will consider two examples of the general set-up. First, in Example~\ref{cor-stable}, we will consider the situation where the noise $W$ is of the stable or tempered stable type. More precisely, let $\{\rho_v\}_{v \in V}$ be given either by 
\begin{align}\label{stabel}
 \rho_v(dx)={}& \big(\1_{\{x\geq 0\}}c_1(v) x^{-\alpha(v)-1} +\1_{\{x< 0\}}c_2(v) \abs{x}^{-\alpha(v)-1}\big)\, dx,
 \shortintertext{or}
 \label{tempered_stabel}
 \rho_v(dx)={}& \rho(dx)=\big(\1_{\{x\geq 0\}}d_1 x^{-\beta-1}e^{-l_1 x} +\1_{\{x< 0\}}d_2 \abs{x}^{-\beta-1} e^{-l_2 \abs{x}}\big)\,dx,
\end{align}
where  $c_1,c_2,\alpha$ are measurable functions from $V$ into $[0,\infty)$, $\alpha(v)\in (0,2)$, $\beta\in (0,2)$, $d_1,d_2\geq 0, d_1+d_2>0$,  $l_1, l_2>0$. Equation~\eqref{stabel} defines the L\'evy measure of  a stable distribution with index $\alpha(v)$ and \eqref{tempered_stabel}  is the L\'evy measure of a tempered stable distribution with a fixed index $\beta$;  see \cite{Cont_Tankov}. 

\begin{example}\label{cor-stable} 
\emph{Suppose that $\{\rho_v\}_{v \in V}$ is given by \eqref{stabel} with $\alpha(v)\in (1+\epsilon,2)$ for some $\epsilon \in (0,1)$, or $\{\rho_v\}_{v \in V}$ is given by \eqref{tempered_stabel} with $\beta\in (1,2)$. Suppose moreover that  $\sigma^2= 0$.   Then, a SIMMA process $X=\pr  X$, given by \eqref{eq:mma}, is of finite variation  if and only if  for $m$-a.e.~$v$, $f(\cdot,v)$ is absolutely continuous  with  derivative $\dot f(\cdot,v)$ satisfying
 \begin{align}\label{}
\begin{dcases} \int_V\int_\R\Big(\frac{c_1(v)+c_2(v)}{2-\alpha(v)}\abs{\dot f(s,v)}^{\alpha(v)}\Big) \, ds \, m(dv)<\infty\qquad \qquad & \text{when $\rho_v$ is given by \eqref{stabel}},\\ \medskip
 \int_V\int_\R\Big(\abs{\dot f(s,v)}^{\beta}\wedge \abs{\dot f(s,v)}^2\Big) \, ds \, m(dv)<\infty\qquad & \text{when $\rho_v$ is given by \eqref{tempered_stabel}}.
\end{dcases}
\end{align}
} 
\end{example}

In the setting of Example~\ref{cor-stable} and $\rho_v$ given by \eqref{stabel} we note that Condition~\eqref{inf-var} of Theorem~\ref{necessary_fv} is satisfied if and only if $\alpha(v)\geq 1$. Moreover, from  \eqref{thy} below, it follows   that Condition~\eqref{eq:u00} of Theorem~\ref{necessary_fv} corresponds to  $\alpha(v)\in (1+\epsilon,2)$. 

\begin{proof}
Let $\{\rho_v\}_{v \in V}$ be given by \eqref{stabel}.  For $v\in V$,  $\int_{-1}^1 \abs{x}\,\rho_v(dx)=\infty$ and hence \eqref{inf-var} is satisfied.  Using \eqref{stabel} a simple calculation shows that
    \begin{align}\label{}
 &\int_\R \big(\abs{xu}\wedge \abs{xu}^2\big)\,\rho_v(dx)=C(v) \abs{u}^{\alpha(v)},\qquad u\in\R,\\
 \shortintertext{where}
& C(v):=\big(c_1(v)+c_2(v)\big)\Big(\frac{1}{\alpha(v)-1}+\frac{1}{2-\alpha(v)}\Big).
\end{align}
A similar calculation shows that   
  \begin{align}\label{thy}
&u \int_{\abs{x}>u} \abs{x}\,\rho_v(dx)=K_0(v)
 \int_{\abs{x}\leq u} x^2\,\rho_v(dx),\qquad u\geq 0,
  \end{align}
 where $ K_0(v)=(2-\alpha(v))/(\alpha(v)-1) $, and  since 
$\alpha(v)\in (1+\epsilon,2)$ by assumption,  \eqref{eq:u00} holds. Hence the result  follows by Theorems~\ref{th:3} and \ref{necessary_fv}.

Assume that  $\rho_v=\rho$ is given by \eqref{tempered_stabel} and note that $\int_{-1}^1 \abs{x}\,\rho(dx)=\infty$.  In the following we will use the notation  $f(u)\sim g(u)$ as $u\to 0$ (or $\infty$), if $f(u)/g(u)\to 1$ as $u\to 0$ (or $\infty$). Moreover, we will use the  asymptotics of the incomplete gamma functions. We have that 
\begin{equation}\label{}
 \rho([-u,u]^c)\sim  (d_1+d_2)\beta^{-1} u^{-\beta} \qquad \text{as } u\to 0, 
\end{equation}
which by Proposition~\ref{remark}  shows that  $\rho$ satisfies \eqref{eq:u00}, keeping in mind that $\int_{\abs{x}>1} x^2\,\rho(dx)<\infty$.     From \eqref{tempered_stabel} we have
 \begin{align}\label{}
 & \int_\R \big(\abs{xu}\wedge \abs{xu}^2\big)\,\rho(dx)\sim \begin{cases} 
 C_1 u^\beta\qquad & \text{as  }u\to \infty,\\ 
 C_2 u^2\qquad & \text{as }u\to 0,\end{cases}
 \end{align}
where $C_1=(d_1+d_2)((\beta-1)^{-1}+(2-\beta)^{-1})$ and $C_2= (d_1l_1^{\beta-2}+d_2l_2^{\beta-2})\Gamma(2-\beta)$, which by Theorems~\ref{th:3}  and \ref{necessary_fv} completes the proof. 
  \end{proof}

 Let $v\mapsto \alpha(v)$ be a measurable function from $V$ into $\R$ and 
 consider   $X=\pr X$ of the form
\begin{equation}\label{def_multi_frac}
 X_t=\int_{\R\times V} \big( (t-s)_+^{\alpha(v)}-(-s)_+^{\alpha(v)}\big) \,W(ds,dv),
\end{equation}
where  $0^0:=0$ and $x_+:=\max\{x,0\}$ for $x\in\R$. We will, as in the rest of this paper, assume that $X$ is well-defined. When $V$ is a one point space, $X$ is called a linear fractional L\'evy process.  Thus, a process $X$ of  form \eqref{def_multi_frac} is a superposition of linear fractional L\'evy processes with (possible) different indexes, and  will therefore be called a  supFLP. 

   \begin{example}\label{cor_frac}
   \textit{Let  $X=\pr X$ be a supFLP of the form \eqref{def_multi_frac}.
   If $\sigma^2= 0$, $\alpha\in [0,\frac{1}{2})$ $m$-a.e.\ and 
  \begin{equation}\label{suf_con_frac}
 \int_V \Big(\int_\R \abs{x}^{\frac{1}{1-\alpha(v)}} \,\rho_v(dx)\Big)\, \big(\tfrac{1}{2}-\alpha(v)\big)^{-1}\,m(dv)<\infty,
\end{equation}
then  $X$ is of finite variation. 
On the other hand, if $X$ is of finite variation, then $m$-a.e., $\sigma^2=0$,  $\alpha\in [0,\frac{1}{2})$  and 
\begin{equation}\label{eq:nece_ex}
 \int_\R \abs{x}^{\frac{1}{1-\alpha(v)}}\,\rho_v(dx)<\infty.\end{equation}
If, in addition,  $\rho$ satisfies \eqref{eq:u00}, then \eqref{suf_con_frac} is satisfied. 
 }
 
 \medskip
To see that  the above example follows from Theorems~\ref{th:3} and \ref{necessary_fv} we need the following general facts about supFLPs $X$ of the form \eqref{def_multi_frac}. Process $X$  is of the form \eqref{eq:mma} with $f(s,v)=f_0(s,v)=s^{\alpha(v)}_+$. Since $X$ is well-defined,  an application of \citet[Theorem~2.7]{Rosinski_spec}  shows that 
 \begin{equation}\label{suf_well}
 \int_V \int_\R \int_\R \Big(\abs{(f(1-s,v)-f_0(-s,v))x}^2\wedge 1\Big)\,ds\,\rho_v(dx)\,m(dv)<\infty.
 \end{equation}
For all $s>0$ there exists $z=z(s,v)\in [s,s+1]$ such that $f(1+s,v)-f(s,v)=\alpha(v)z^{\alpha(v)-1}$. By \eqref{suf_well}  it follows  that $\alpha <\frac{1}{2}$ $m$-a.e.\ and since  $z^{\alpha-1}\geq (s+1)^{\alpha-1}$, \eqref{suf_well} shows that   
 \begin{equation}
 \label{eq:al} 
 \int_V \int_{\abs{x\alpha(v)}>1}  \Big(\frac{\abs{\alpha(v)x}^{\frac{1}{1-\alpha(v)}}}{1-2\alpha(v)}\Big)\,\rho_v(dx)\,m(dv)<\infty,
\end{equation} 
which implies that  
\begin{equation}\label{int-m-frac}
 \int_{\abs{x}>1} \abs{x}^{\frac{1}{1-\alpha(v)}}\,\rho_v(dx)<\infty\quad \text{for }m\text{-a.e.\ } v. 
\end{equation} 
For $v\in V$,  $f(\cdot,v)$ is absolutely continuous if and only if $\alpha(v)>0$ and in this case $\dot f(s,v)=\alpha(v)s^{\alpha(v)-1}_+$.  For $\alpha(v)\in (0,\frac{1}{2})$, a simple calculation shows that 
\begin{align}\label{sim_cal}
  \int_\R \big(\abs{\dot f(s,v)x}^2\wedge \abs{\dot f(s,v)x}\big) \,ds= 
  \abs{x}^{\frac{1}{1-\alpha(v)}}\Big[ \abs{\alpha(v)}^{\frac{1}{1-\alpha(v)}}\Big( \frac{1}{\alpha(v)}+\frac{1}{1-2\alpha(v)}\Big)\Big].
    \end{align} 
    The square bracket in \eqref{sim_cal} is, for  $\alpha(v)\in (0,\frac{1}{2})$, bounded from above and below by two constants $c_1,c_2>0$  times  $(\tfrac{1}{2}-\alpha(v))^{-1}$, which shows that 
    \begin{equation}\label{sim_cal1}
\frac{ c_1 \abs{x}^{\frac{1}{1-\alpha(v)}}}{\tfrac{1}{2}-\alpha(v)}\leq \int_\R \big(\abs{\dot f(s,v)x}^2\wedge \abs{\dot f(s,v)x}\big) \,ds\leq \frac{ c_2 \abs{x}^{\frac{1}{1-\alpha(v)}}}{\tfrac{1}{2}-\alpha(v)}.
\end{equation}

\begin{proof}[Proof of Example~\ref{cor_frac}]
Let $f(s,v)=f_0(s,v)=s^{\alpha(v)}_+$. We may and do consider the following two cases separately: $\alpha(v)=0$ for all $v\in V$, and $\alpha(v)\neq 0$ for all $v\in V$.   If  $\alpha(v)=0$ for all $v\in V$, then, $X_t=W((0,t]\times V)$  is a L\'evy process with L\'evy measure $\nu(dx)=\,\rho_v(dx)\,m(dv)$ and Gaussian component $\int_V \sigma^2(v)\,m(dv)$. Hence $X$ is  of  finite variation if and only if $\int_{V}\int_\R \abs{x}\,\rho_v(dx)\,m(dv)<\infty$ and $\sigma^2=0$ $m$-a.e., cf.\ \cite[Theorem~21.9]{Sato}. Thus, in what follows  we will assume that  $\alpha(v)\neq 0$ for all $v\in V$. 
 
 Assume that $\alpha\in (0,\frac{1}{2}),\ \sigma^2=0$ $m$-a.e.\ and \eqref{suf_con_frac} is satisfied. For $m$-a.e.\ $v$,  $f(\cdot,v)$ is absolutely continuous and by \eqref{sim_cal1},  $\dot f(\cdot,v)$ satisfies \eqref{eq:Cf}, which  by  Theorem~\ref{th:3} shows that   $X$ is of finite variation. 
 
 On the other hand, assume that $X$ is of finite variation.  By a symmetrization argument we may consider the cases where $W$ is centered Gaussian or has no Gaussian component separately.  In the Gaussian case we have $\sigma^2>0$ $m$-a.e.\  by \eqref{eq:pur-sto}, and therefore  \eqref{inf-var} holds.   For $m$-a.e.\ $v$, $s\mapsto f(s,v)$ is absolutely continuous with a derivative $\dot f(s,v)=\alpha(v)s^{\alpha(v)-1}_+$ satisfying \eqref{gen_eq}, cf.\    Theorem~\ref{necessary_fv}.  Hence    $\alpha(v)>0$ and   by \eqref{sim_cal1} 
 \begin{equation}\label{eq_er2}
  \int_V \Big(\int_0^\infty  \abs{s^{\alpha(v)-1}}^2\,ds\Big)\,\abs{\alpha(v)}^2\sigma^2(v)\,m(dv)<\infty, 
 \end{equation}
 which implies that $\sigma^2=0$ $m$-a.e.
 In the purely non-Gaussian case, \citet[Theorem~4]{Rosinski_sum_rep} shows that  $f(\cdot,v)$ is of finite variation for $m$-a.e.\ $v$. Hence $\alpha\geq 0$ and by assumption $\alpha>0$. 
  Thus for $m$-a.e.\ $v$, $f(\cdot,v)$ is  absolutely continuous and  by \eqref{sim_cal1} and the below Remark~\ref{remark-con}  we have
   \begin{align}\label{eq_er}
  \int_\R \Big(\frac{\abs{x}^{\frac{1}{1-\alpha(v)}}}{1\vee x^2}\,\Big)\,\rho_v(dx)<\infty\qquad \text{for $m$-a.e.\ }v,
\end{align}
which combined with \eqref{int-m-frac} shows \eqref{eq:nece_ex}.  Finally, if  $\rho$ satisfies \eqref{eq:u00} then Remark~\ref{remark-con} and \eqref{sim_cal1} show that \eqref{suf_con_frac} is satisfied. This   completes the proof. 
   \end{proof}
 \end{example}

In the special case where  $V$ is a one point space, i.e.\ $X$ is a fractional L\'evy process,  Example~\ref{cor_frac} shows that $X$ is of finite variation  if and only if $\sigma^2=0, \alpha\in [0,\frac{1}{2})$ and \eqref{eq:nece_ex} is satisfied.  This completes  \cite[Corollary~5.4]{Basse_Pedersen}  and parts of \cite[Theorem~2.1]{Be_Li_Sc}.

\bigskip

%%%%%%%%%%%%%%%%%%%
\section{Proofs of Theorems~\ref{th:3} and \ref{necessary_fv}}
\label{proof_fv}
%%%%%%%%%%%%%%%%%%%%

We will start by showing Theorem~\ref{th:3}.

\begin{proof}[Proof of Theorem~\ref{th:3}]
 Let $B=\{v:\dot f(\cdot,v) = 0\ \lambda\text{-a.e.}  \}$. By  \eqref{eq:Cf}, $\int_{\abs{x}>1} \abs{x}\,\rho_{\cdot}(dx)<\infty$ $m$-a.e.\ on $B^c$, and since $f(\cdot,v)$ is constant  for $v\in B$ we may and do assume that  $\int_{\abs{x}>1} \abs{x}\,\rho_{\cdot}(dx)<\infty$ $m$-a.e. This allows us to write $W$ as $W=W_0+\mu$, where $W_0$ is a centered random measure and 
 $\mu$ is a deterministic measure. 
  To show that $\pr X$ has absolutely continuous sample paths, define a measurable process $\pr{Y^0}$ by
\begin{equation} \label{def_y0}
Y_t^0= \int_{\R\times V} {\dot f}(t-s,v) \, W_0(ds, dv).
\end{equation}
By the assumptions \eqref{gen_eq} and \eqref{eq:Cf},   a stochastic Fubini theorem,  see \cite[Remark~3.2]{QOU}, shows that  process $Y^0$ is well-defined and  for all $a<b$,
\begin{align*}
 \int_a^b Y_t^0 \, dt =  {}&\int_{\R\times V}\Big( \int_a^b  {\dot f}(t-s,v) \, dt \Big)\, W_0(ds, dv) \\   = {}&  \int_{\R\times V} \big(f(b-s,v)-f(a-s,v)\big)\,W_0(ds,dv)
\end{align*}
with all integrals well-defined. 
By linearity,
\begin{equation}\label{}
h(t) :=  \int_{\R\times V} \big(f(t-s,v)-f(-s,v)\big)\,\mu(ds,dv),\qquad t\in\R,
\end{equation}
is well-defined as well.   Using that $h(t)=h(t+u)-h(u)$ for all $u,t\in\R$ and that $h$ is measurable,  a standard argument shows that 
 $h(t)=t  h(1)$.  Thus, with $Y_t:=h(1)+Y^0_t$, we have with  probability 1, 
\begin{equation} \label{}
X_t= X_0+ \int_0^t Y_u \, du, \qquad t \in \R,
\end{equation}
which proves Theorem~\ref{th:3}. 
\end{proof}

\begin{proof}[Proof of Corollary~\ref{cor_var}]
  Corollary~\ref{cor_var} follows by the estimates given in \citet{Rosinski_L1_norm}, Corollary~1, used on $Y^0_t$ in \eqref{def_y0}. 
 \end{proof}

To prove  Theorem~\ref{necessary_fv} we need the following 
Lemmas~\ref{int_ID} and \ref{co_type_2} about general symmetric infinitely divisible processes.  
Let $T$ denote a countable set and $X=\xT$ be a symmetric infinitely divisible process without Gaussian component. Let $\R^T$ be equipped with the product topology,  $\R^{(T)}$ denote its 
the topological dual space,  and $\inner{\cdot,\cdot}$ be the canonical bilinear form on $\R^{(T)}\times \R^T$.  For each $y\in\R^{(T)}$ there exist $n\in \N$, 
$(\alpha_i)_{i=1}^n\subseteq \R$ and $(t_i)_{i=1}^n\subseteq T$ such that 
$\inner{y,x}=\sum_{i=1}^n \alpha_i x_{t_i}$ for all $x\in \R^T$.
Let $\nu$ be the L\'evy measure of $X$, that is, $\nu$ is a symmetric Borel measure on $\R^T$ with $\nu(\{0\})=0$ and 
$\int (1\wedge x(t)^2)\,\nu(dx)<\infty$ for all $t\in T$ such that for all $y\in \R^{(T)}$,
\begin{align}\label{LK}
\E e^{i\inner{y,X}}=\exp\Big(\int_{\R^T} \big(\cos(\inner{y,x})-1\big)\,\nu(dx)\Big).
\end{align}

Let $\morf{h}{[0,\infty)}{[0,\infty)}$ be a submultiplicative function, i.e., there exists a constant $c>0$ such that 
\begin{equation}\label{sub-multi}
h(x+y)\leq c h(x)h(y),\quad x,y\geq 0.
\end{equation} 
Assume, moreover, that $h$ is increasing,  and for all $\epsilon >0$ there exists $a_\epsilon>0$ 
such that  $h(x)\leq a_\epsilon e^{\epsilon x}$ for all $x\geq 0$. Let $h(\infty)=\infty$. 
The key example is $\map{h}{x}{(x\vee 1)^p}$ for  $p>0$, where $x\vee 1=\max\{x,1\}$. If $q$ is a lower semicontinuous pseudonorm on $\R^T$ such that $q(X)<\infty$ a.s., Lemma~2.1 in \cite{Ros_Sam} shows that 
there exists an $r_0\in (0,\infty)$ such that $\nu(x\in\R^T:q(x)\geq r_0)<\infty$. 

\begin{lemma}\label{int_ID}
Let $T$ be a countable set, $X=\xT$ be a symmetric  infinitely divisible process of the form \eqref{LK} and  $\morf{q}{\R^T}{[0,\infty]}$ be
a lower-semicontinuous pseudonorm such that $q(X)<\infty$ a.s. For all $r_0>0$ such that $\nu(x\in \R^T:q(x)\geq r_0)<\infty$ we have
\begin{equation}\label{int_iff}
 \int_{\{q(x)\geq r_0\}} h(q(x))\,\nu(d x)<\infty\ \ \text{if and only if}\    \  \E h(q(X))<\infty. 
\end{equation}
\end{lemma} 

Lemma~\ref{int_ID} in the finite dimensional case, i.e., $\rm{Card}(T) < \infty$, follows from  \citet[Theorem~25.3]{Sato}.    The case $\rm{Card}(T) = \infty$ requires some minor changes. For example, we use \citet[Lemma~2.2]{Ros_Sam} instead  of  the Lemmas~25.6 and 25.7 in \cite{Sato}, since the latter do not extend to an infinite dimensional case. 

\begin{proof}[Proof of Lemma~\ref{int_ID}]
%We claim that 
%\begin{equation}\label{nu-nul}
%\nu(x\in\R^T:q(x)=\infty)=0.
%\end{equation}
%To show \eqref{nu-nul} let $\tilde X$ denote an independent copy
%of $X$ and define $\tilde X=X- X'$. Then $\tilde X$ is a symmetric infinitely 
%divisible process with L\'evy measure $\tilde \nu(A)=\nu(A)+\nu(-A)$. Let $\tilde \nu_1=\tilde \nu_{ |\{q\leq  r_0\}}$, 
% $\tilde \nu _2=\tilde \nu_{|\{q>r_0\}}$, and  $\tilde X^1$ and 
%$\tilde X^2$ two independent and symmetric infinitely divisible processes with L\'evy measures $\tilde \nu _1$ and $\tilde \nu _2$. By symmetry
%\begin{equation}
%\tilde X\dist \tilde X^1+\tilde X^2\dist \tilde X^1-\tilde X^2
%\end{equation}
%where $\dist$ means equality in finite dimensionally distributions.  Since 
%$q(\tilde X)<\infty$ a.s.\ we have that 
%\begin{equation}
%q(\tilde X^2)\leq q\Big(\frac{\tilde X^2+\tilde X^1}{2}\Big)+q\Big(\frac{\tilde X^2-\tilde X^1}{2}\Big)<\infty,\qquad \text{a.s.}
%\end{equation}
%For each measure $m$ on $\R^T$  and $k\in \N$ let $m^{\otimes k}$ denote the $k$-fold convolution of $m$ and   $m^{\otimes 0}:=\delta_{0}$.
% Since 
%\begin{equation}
%\P(\tilde X^2\in A)=e^{-\tilde \nu_2(\R^T)}\sum_{k=0}^\infty \frac{\tilde \nu_2^{\otimes k}(A)}{k!},\qquad A\in \B(\R^T),
%\end{equation} 
%we have that $0=\P(q(\tilde X^2)=\infty)\geq e^{-\tilde \nu_2(\R^T)}\tilde \nu_2(q=\infty)$. That is,   $0=\nu_2(q=\infty)$ which shows  \eqref{nu-nul}. 
% 
Let $\nu_1:=\nu_{|\{q< r_0\}}$, $\nu_2:=\nu_{|\{q\geq r_0\}}$, and  
 $X^1$ and $X^2$ be two independent symmetric infinitely divisible processes  such that for all $\beta\in \R^{(T)}$, 
\begin{align}
\E[e^{i\inner{\beta,X^1}}]={}&\exp\Big(
\int_{\R^T} \big(\cos(\inner{y,x})-1\big)\,\nu_1(dx)\Big),\\
\E[e^{i\inner{\beta,X^2}}]={}&\exp\Big(
\int_{\R^T} \big(\cos(\inner{y,x})-1\big)\,\nu_2(dx)\Big).
\end{align} 
By convexity, 
\begin{align}
q(X^1)\leq \frac{1}{2}\Big(q(X^1+X^2)+ q(X^1-X^2)\Big),
\end{align}
which shows that  $q(X^1)<\infty$ a.s.\ due to the fact that $X\dist X^1+X^2\dist X^1-X^2$ where $\dist$ denotes equality of finite dimensionally distributions. 

To show the \emph{only if}-implication assume that the left-hand side of \eqref{int_iff} is satisfied.
Since  $q(X^1)<\infty$ a.s.\ and $\nu_1(q\geq  r_0)=\nu(\emptyset)=0$, Lemma~2.2 in \cite{Ros_Sam} shows that there 
exists an $\epsilon>0$ such that $\E e^{\epsilon q(X^1)}<\infty$. By assumption there exists $a_\epsilon>0$ such that $h(x)\leq a_\epsilon e^{\epsilon x}$ for all $x\geq 0$  and hence 
\begin{align}
\E[h(q(X^1))]\leq a_\epsilon \E[e^{\epsilon q(X^1)}]<\infty.
\end{align} 
For any $k\in\N$ let $\nu^{\otimes k}_2$ denote the $k$-fold convolution of $\nu_2$ and $\nu_2^{\otimes 0}:=\delta_0$. We may and do assume that the constant $c$ from \eqref{sub-multi} satisfies $c\geq 1$ and hence   
\begin{align}
\E[h(q(X^2))]={}& e^{-\nu_2(\R^T)}\sum_{k=0}^\infty \frac{\int h(q(x))\,\nu_2^{\otimes k}(dx)}{k!} \\  \leq {}&  e^{-\nu_2(\R^T)}h(0)+\sum_{k=1}^\infty \frac{c^{k-1}}{k!}\Big(\int
h(q(x))\,\nu_2(dx)\Big)^k \\ \leq {}& e^{-\nu_2(\R^T)}h(0)+\exp\Big(c\int_{\{q\geq
  r_0\}} h(q(x))\,\nu(dx)\Big)<\infty. 
\end{align}
Since $h$ is submultiplicative and increasing, 
\begin{align}\label{sub_add}
\E[h(q(X))]=\E[h(q(X^1+X^2))]\leq c \E[h(q(X^1))]\E[h(q(X^2))]<\infty
\end{align} 
which shows  that the right-hand side of \eqref{int_iff} is satisfied. 

To show the \emph{if}-implication assume that $\E [h(q(X))]<\infty$. 
Since  $q(x)<\infty$ for $\P_{X^1}$-a.a.\ $x$ and 
\begin{equation}\label{}
\infty>\E[h(q(X))]=\int_{\R^T} \E[h(q(x+X^2))]\,\P_{X^1}(dx),
\end{equation} 
there exists $x\in \R^T$ with $q(x)<\infty$ such that 
$\E[ h(q(x+X^2))]<\infty$. Hence 
\begin{equation}
\E[h(q(X^2))]\leq c \E[h(q(x+X^2))] h(q(x))<\infty,
\end{equation} 
and the left-hand side of \eqref{int_iff} follows by the inequality
\begin{align}
\E[h(q(X^2))]=e^{-\nu_2(\R^T)}\sum_{k=0}^\infty \frac{\int h(q(x))\,\nu_2^{\otimes k}(dx)}{k!}\geq e^{-\nu_2(\R^T)}\int_{\{q(x)\geq r_0\}} h(q(x))\,\nu(dx).
\end{align}
\end{proof}

\begin{lemma}\label{co_type_2}
  Let   $N\in \N$,  $T=\{k2^{-n}:n\in \N,\, k=0,\dots,N2^n\}$ and for $\morf{f}{T}{\R}$ define 
  \begin{equation}\label{}
  \norm{f}_{BV[T]}=\sup_{n\in\N}
  \sum_{k=1}^{N2^n}\abs[\big]{f(k2^{-n})-f((k-1)2^{-n})}.
\end{equation} 
For any  infinitely divisible process $X=(X_t)_{t\in T}$ of the form \eqref{LK} with $\norm{X}_{BV[T]}<\infty$ a.s.\  we have 
  \begin{equation}\label{}
 \int_{\R^T} \big(1\wedge \norm{x}^2_{BV[T]} \big)\,\nu(dx)<\infty. 
\end{equation}
 \end{lemma}
 
  It can be shown that $BV[T]$ is a Banach space of cotype 2, however, it is not separable so \citet[Theorem~2.2]{cotype2} does not apply to this situation.   To prove Lemma~\ref{co_type_2}  we  use  \citet[Theorem~4.9]{Rosinski_spec}  and  \citet[Proposition~2]{Rosinski_sum_rep}  which gives a series representation of  $X$. Using the series representation,  Lemma~\ref{co_type_2}  follows  along the lines of  Proposition~5.6 in \citet{Basse_Pedersen}.

We are now ready to prove Theorem~\ref{necessary_fv}.

\begin{proof}[Proof of Theorem~\ref{necessary_fv}]
We need to show the following three cases \emph{(a)}: \eqref{gen_eq} and  \eqref{eq:1suf}  hold under no additional restrictions on $\rho$,    \emph{(b)}:  \eqref{fdot_int} holds under \eqref{eq:u0}, \emph{(c)}:  \eqref{eq:Cf} holds under  \eqref{eq:u00}. We will start by showing \emph{(b)} and at the end of the proof   deduce \emph{(a)} and \emph{(c)} from it. 

\emph{(b)}: Assume  \eqref{eq:u0}.  Using the monotonicity
\begin{align}\label{}
   \sum_{k=1}^{2^n} & |X_{k2^{-n}} - X_{(k-1)2^{-n}}| \\ 
  & \le \sum_{k=1}^{2^n}  \big( \abs{X_{(2k)2^{-n-1}} - X_{(2k-1)2^{-n-1}}}+\abs{X_{(2k-1)2^{-n-1}} - X_{(2k-2)2^{-n-1}}} \big) \\
  &  = \sum_{k=1}^{2^{n+1}} |X_{k2^{-n-1}} - X_{(k-1)2^{-n-1}}|
\end{align}
and the stationarity of increments of $X$, we get by \eqref{eq:fvX} 
\begin{align} 
\E \|X\|_{BV[0,1]} &= \E \bigg[\lim_{n\to \infty}  \sum_{k=1}^{2^n}  \abs{X_{k2^{-n}} - X_{(k-1)2^{-n}}}\bigg] = \lim_{n\to \infty}  \sum_{k=1}^{2^n}  \E \abs{X_{k2^{-n}} - X_{(k-1)2^{-n}}} \\
&= \sup_{n\in \N} \E\abs{2^n(X_{2^{-n}} - X_{0})}. \label{eq:mBV}
\end{align}

%Let $W'$ be an independent copy of $W$ and set $\bar{W}=W-W'$. Moreover, let $\bar{X}$ be a separable process defined by \eqref{eq:mma} with $W$ replaced by $\bar{W}$. 
%Then $\bar{X}$ has a.s.\ sample paths of finite variation. By considering $\bar{X}$ instead of $X$ we may and will assume in the following that 
By a symmetrization argument, see Section~\ref{sec_MA},  we may and will assume that $W$ is a symmetric random measure and,  in particular,  we may and will assume that  $\theta=0$ and $\rho_v$ are symmetric L\'evy measures. Decompose $W$ as $W=W_P+W_G$ where $W_P$ and $W_G$ are random measures of the form \eqref{eq:M} with $(\theta_P,\sigma^2_P,\rho_P)=(0,0,\rho)$ and  $(\theta_G,\sigma^2_G,\rho_G)=(0,\sigma^2,0)$, respectively. Let $X^P=(X_t^P)_{t\in \R}$ and $X^G=(X_t^G)_{t\in \R}$ be  processes of the form \eqref{eq:mma} with $W$ replaced by $W_P$ and  $W_G$, respectively. Processes $X^P$ and $X^G$ are chosen separable with the dyadics $\D$ as their separant. By symmetry, $X^P$ and $X^G$ have sample paths  of finite variation almost surely.  We will divide the proof into  the following  Steps~\emph{(b1)}--\emph{(b3)}.  In Steps~\emph{(b1)} and \emph{(b2)} we will consider, respectively,  the non-Gaussian and Gaussian case separately, and in the Step~\emph{(b3)} we will  deduce the general case  from Steps~\emph{(b1)}--\emph{(b2)}. 

Step~\emph{(b1)}: For any $v\in V$ let 
\begin{equation} \label{}
\xi_v(u) := \int_{\R} ( |ux|^2 \wedge |ux|) \, \rho_v(dx),\quad u\in \R.
\end{equation}
Then $\xi_v$ is symmetric,  increasing, and comparable with a convex function $\tilde{\xi_v}$ given by
\begin{equation} \label{}
\tilde{\xi_v}(u) = \int_{\R} (|ux|^2 \1_{\{|ux|\le 1\}} + (2|ux|-1)\1_{\{|ux|>1\}}) \, \rho_v(dx).
\end{equation}
Indeed,  $\tilde{\xi_v}(u)/2 \le \xi_v(u) \le \tilde{\xi_v}(u)$, $u\ge 0$.
By Corollary~1.1 in \citet{Rosinski_L1_norm}
\begin{equation} \label{eq:L1}
\frac{1}{4} \min\{I_n, I_n^{1/2}\} \le \E \abs{2^n(X_{2^{-n}}^P - X_{0}^P)}\le \frac{5}{4} \max\{I_n, I_n^{1/2}\},
\end{equation}
where 
\begin{equation} \label{def_fn7}
I_n = \int_{\R}\int_V \xi_v(f_n(s,v)) \, m(dv)\, ds\quad \text{and}\quad 
f_n(s,v)=2^n[f(2^{-n}-s, v) - f(-s, v)]. 
\end{equation}
(For symmetric infinitely divisible random variables \cite[Corollary~1]{Rosinski_L1_norm} remains true without the first moment condition.) 
In view of \eqref{eq:mBV} and \eqref{eq:L1}, 
\begin{equation}\label{int_var_iff}
  \E \|X^P\|_{BV[0,1]} < \infty\ \text{ if and only if } \  \sup_{n\in\N} I_n < \infty. 
\end{equation}
For $N\in \N$ set  $T=\{k2^{-n}:n\in \N,\, k=0,\dots,N2^n\}$. Then  $(X_t^P)_{t\in T}$  is an infinitely divisible process of the form \eqref{LK} with L\'evy measure $\nu$ determined by  
\begin{equation}
\nu(p_{t_1,\dots,t_n}^{-1}(A))=\int_\R \int_V\int_\R  \1_A(xf(t_1-s,v),\dots,xf(t_n-s,v))\,\rho_v(dx)\,m(dv)\,ds
\end{equation}
for all $n\in \N$, $t_1,\dots, t_n\in T$, $A\in \B(\R^n)$ and with $\morf{p_{t_1,\dots,t_n}}{\R^T}{\R^n}$   given by $x\mapsto (x(t_1),\dots, x(t_n))$, see \cite[Theorem~2.7]{Rosinski_spec}. 
Since  the dyadic numbers are a separant for $(f(t-\cdot,\cdot))_{t\in\R}$, we have by  Lemma~\ref{co_type_2} that 
\begin{equation} \label{eq:co2}
\int_{\R} \int_V \int_{\R} \left(1 \wedge \|xf(\cdot - s,v)\|_{BV[0,N]}^2 \right) \, \rho_v(dx) \, m(dv)\, ds < \infty.
\end{equation}
Since \eqref{eq:co2} holds for all $N\in \N$  there exists a measurable set  $V_0\in  \V$ with $m(V\setminus V_0)=0$ such that for every $v\in V_0$ and $t>0$
\begin{equation} \label{eq:co2.1}
\int_{\R} \left(1 \wedge \|f(\cdot - s,v)\|_{BV[0,t]}^2 \right) \,  ds < \infty.\end{equation}
We will show that 
\begin{equation}\label{newdefk}
 k^*(v):=\sup_{s\in\R} \norm{f(\cdot-s,v)}_{BV[0,1]}<\infty,\qquad  v\in V_0. 
\end{equation}
To do this notice that 
\begin{equation} \label{}
\|f(\cdot - s,v)\|_{BV[0,t]} = \|f(\cdot,v)\|_{BV[-s,t-s]} = k(t-s, v) - k(-s, v),
\end{equation}
where 
\begin{equation} \label{}
k(u,v) = \begin{cases}
 \|f(\cdot,v)\|_{BV[0,u]}   & \text{if } u\ge 0, \\ 
  -\|f(\cdot,v)\|_{BV[u,0]}   & \text{if } u< 0. 
\end{cases}
\end{equation}
For each $v\in V_0$, $u \mapsto k(u,v)$ is a nondecreasing function.
To show \eqref{newdefk}  fix $v\in V_0$ and let us for the moment suppress $v$. Let $h(s)=\abs{k(1-s)-k(-s)}$.  For contradiction assume that $h$ is 
unbounded. Since $h$ is locally bounded there exists a sequence $(a_n)_{n\in\N}$ converging to either
 $\infty$ or $-\infty$ (say, $\infty$) such that
 $h(a_n)\geq 1$ for all $n\in\N$. By passing to a subsequence we
 may   assume that $a_n+1\leq a_{n+1}$ for all $n\in\N$. For $s\in [a_n,a_n+1]$ we have 
 \begin{equation}
 k(2-s)-k(-s)\geq k(1-a_n)-k(-a_n)=h(a_n)\geq 1.
 \end{equation} 
 Thus, 
 \begin{equation}
 \int_\R \big(1\wedge [k(2-s)-k(-s)]^2\big)\,ds\geq \sum_{n=1}^\infty \int_{a_n}^{a_n+1} \big(1\wedge [k(2-s)-k(-s)]^2\big)\,ds
 \geq \sum_{n=1}^\infty 1=\infty, 
 \end{equation} which  contradicts \eqref{eq:co2.1} and completes the proof of \eqref{newdefk}. 
 
 By \eqref{eq:u0} there  exist two measurable functions $\morf{u_0}{V}{[0,\infty)}$ and $\morf{K_0}{V}{(0,\infty)}$ such that for $m$-a.e.\ $v$
 \begin{equation}\label{eq:u0-1}
u\int_{\abs{x}>u} \abs{x}\,\rho_v(dx)\leq K_0(v)\int_{|x|\le u} x^2 \, \rho_v(dx)\quad \text{for  }u\geq u_0(v),
\end{equation}
which implies that 
 \begin{align}\label{eq:rho1}
{}&  \int_{\abs{ux}>1} \abs{xu}\,\rho_v(dx)\leq K_0(v)\int_{\R} (\abs{xu}^2\wedge 1)\,\rho_v(dx)\qquad \text{for }\abs{u}\leq  1/u_0(v).
\end{align}
     For arbitrary but fixed $k\in \N$ define 
\begin{equation} \label{}
V_k = \{v \in V_0: k^*(v) \le k,\ K_0(v)\leq k, u_0(v)\leq k\},
\end{equation}
and let $\pr {X^k}$ be given by 
\begin{equation}\label{}
 X^k_t= \int_{\R\times V_k} \big(f(t-s, v)- f_0(-s, v)\big)\, W_P(ds, dv),\qquad t\in\R.
\end{equation}
By a symmetrization argument, $\|X^k\|_{BV[0,1]} < \infty$ a.s. 
We will show that
\begin{equation} \label{eq:mN}
\E \|X^k\|_{BV[0,1]} < \infty.
\end{equation}
To this end it is enough, according to Lemma~\ref{int_ID},  to prove that 
\begin{equation} \label{int_eq}
\int_{\{(s,v,x) \in \R\times V_k\times \R: \, \|xf(\cdot - s,v)\|_{BV[0,1]}>k^2\}}  \|xf(\cdot - s,v)\|_{BV[0,1]}\, ds \, \rho_v(dx)\,m(dv)  < \infty.
\end{equation}
 For  $v\in V_k$,  
\begin{equation}\label{}
  \norm{f(\cdot-s,v)}_{BV[0,1]} k^{-2}\leq k^*(v) k^{-2}\leq k^{-1}\leq 1/u_0(v).
\end{equation}
Thus applying \eqref{eq:rho1} on $u= \norm{f(\cdot-s,v)}_{BV[0,1]} k^{-2}$ shows  that   the left-hand side of \eqref{int_eq} is less than or equal to
\begin{align}\label{less-inf}
 k^3\int_{ \R}
 \int_{V_k}\int_\R\Big(1\wedge \|xf(\cdot - s,v)\|_{BV[0,1]}^2\Big)\, \rho_v(dx)\, m(dv) \,ds
\end{align}
 which is finite  by \eqref{eq:co2}. This completes the proof of \eqref{eq:mN}.
 
Since $\E \|X^k\|_{BV[0,1]} < \infty$, \eqref{int_var_iff} shows that 
\begin{equation}\label{xi_int}
 \sup_{n\in\N} \int_{V_k}\int_{\R} \xi_v(f_n(s,v))\, ds\, m(dv)<\infty,
\end{equation} 
where $f_n$ are given by \eqref{def_fn7}. Set 
\begin{equation}\label{def-J}
  J=\big\{v\in V: \int_{-1}^1 \abs{x}\,\rho_v(dx)=\infty\big\},
  \end{equation}
and choose $(A_k)_{k\in\N}\subseteq \V$ such that $A_k\uparrow V$ and  $m(A_k)<\infty$ for all $k\in\N$. Let $\lambda_k:=\lambda_{| [-k,k]}$ and $m_k:=m_{| A_k\cap V_k\cap J}$. Note that $\lambda_k\otimes m_k$ is a finite measure.  
For $v\in J$, we have by the monotone convergence theorem that 
\begin{equation}\label{}
  \frac{\xi_v(x)}{x}=\int_\R \big(\abs{u^2 x}\wedge \abs{u}\big)\rho_v(du)\nearrow \int_\R \abs{u}\,\rho_v(du)=\infty\qquad \text{as }x\nearrow \infty.
\end{equation}
Hence  for all $k\in\N$ there exists,  by Egorov's Theorem (see \cite{Strange_func}, Chapter~9, Theorem~1), $B_k\in  \V$ with $m_k(B^c_k)<1/k$  such that   for all $C>0$ there exists $K>0$ such that for all $v\in B_k$, $\inf_{x>K} (\xi_v(x)/x)\geq C$. With  $\tilde m_k:=m_{k | B_k}$ we have that  
\begin{align}\label{}
 \int_{\abs{f_n}>K} \abs{f_n}\, d(\lambda_k\otimes \tilde m_k)={}&\int_{\abs{f_n}>K} \xi_v(f_n(s,v))\frac{\abs{f_n(s,v)}}{\xi_v(f_n(s,v))}\, (\lambda_k\otimes \tilde m_k)(ds,dv) \\  \leq {}& \int_{\abs{f_n}>K} \xi_v(f_n(s,v))\big(\sup_{x>K}\frac{x}{\xi_v(x)}\big)\,  (\lambda_k\otimes \tilde m_k)(ds,dv) \\  \leq {}& \frac{1}{C} \int \xi_v(f_n(s,v))\,  (\lambda_k\otimes \tilde m_k)(ds,dv),
\end{align}
which shows that
\begin{equation}\label{xi-ui3}
 \sup_{n\in\N}  \int_{\abs{f_n}>K} \abs{f_n}\, d(\lambda_k\otimes \tilde m_k)\leq \frac{1}{C}\sup_{n\in\N} \int_{V_k} \int_{\R}\xi_v(f_n(s,v))\, ds\,  m(dv).
\end{equation}
By \eqref{xi_int}--\eqref{xi-ui3} we conclude that  $\{f_n:n\in\N\}$ is uniformly integrable with respect to $\lambda_k\otimes \tilde m_k$. 
Therefore,  by 
the Dunford-Pettis Theorem, see \cite{Dunford}, IV.8, Corollary~11,  there exists a subsequence $(n_j)_{j\in\N}$ and a $h\in L^1(\lambda_k\otimes \tilde m_k)$ such that $\lim_j f_{n_j}=h$ in $\sigma(L^1,L^\infty)$. For all $A\in \V$ with $A\subseteq A_k\cap V_k\cap B_k
\cap J$ and for $(\lambda\otimes \lambda)$-a.e.\ $(s,t)$ with $-k \leq s<t\leq k$, 
\begin{align}\label{eq_h}
{}& \int_A \Big(\int_s^t  h(u,v)\,du\Big) \, m(dv)=\lim_{j\to\infty} \int_A \Big(\int_s^t  f_{n_j}(u,v)\,du\Big) \, m(dv)
\\{}& \quad =   \lim_{j\to\infty} 2^{n_j}\Big[\int_{s+2^{n_j}}^{t+2^{n_j}}  \Big(\int_A f(u,v)\, m(dv)\Big) \,du-\int_{s}^{t}  \Big(\int_A f(u,v)\, m(dv)\Big) \,du\Big]\\ {}& \quad 
=\lim_{j\to\infty} 2^{n_j}\int_{t}^{t+2^{n_j}}  \Big(\int_A f(u,v) \, m(dv)\Big)\,du-\lim_{j\to\infty}2^{n_j}\int_{s}^{s+2^{n_j}}  \Big(\int_A f(u,v)\, m(dv)\Big) \,du\\ \label{eq_f}  {}& \quad  = \int_{A} \big(f(t,v)-f(s,v)\big) \,m(dv).
\end{align}
Since  \eqref{eq_h}--\eqref{eq_f}  is satisfied for $k$ arbitrary we have for  $(\lambda\otimes\lambda\otimes m)\text{-a.e.\ } (t,s,v)\in\R\times\R\times J$,
\begin{equation}\label{}
 f(t,v)-f(s,v)=\int_s^t  h(u,v)\, du,
\end{equation}
which shows that for $m$-a.e.\ $v\in J$, $f(\cdot,v)$ is absolutely continuous with derivative $h(\cdot,v)$. 

Step~\emph{(b2)}:  Set  
\begin{equation}\label{def-G}
G=\{v\in V:\sigma^2(v)>0\}.
\end{equation}
   By Gaussianity,  \citet{Fernique_int}  shows that $\E \norm{X^G}_{BV[0,1]}<\infty$. Let $f_n$ be given by 
\eqref{def_fn7}.
   As in \eqref{eq:mBV} we have that 
  \begin{align}\label{bound-f-n1}
 \E  \norm{X^G}_{BV[0,1]}={}& \sup_{n\in\N} \big(2^n \E \abs{X_{1/2^n}^G-X_0^G}\big)
= \sqrt{\frac{2}{\pi}} \sup_{n\in\N} \big(2^n \norm{X_{1/2^n}^G-X_0^G}_{L^2}\big)
 \\  \label{bound-f-n2}  ={}&\sqrt{\frac{2}{\pi}} \Big(
  \sup_{n\in\N} \int_V\int_\R | f_n(s,v)|^2\sigma^2(v)\, ds\,m(dv)\Big)^{1/2},
\end{align} 
where in the second equality we have used the identity $\norm{U}_{L^1}=(2/\pi)^{1/2} \norm{U}_{L^2}$ for  centered Gaussian random variables $U$. 
Let $\mu(ds,dv)=ds\,\sigma^2(v)\, m(dv)$ be a measure on $\R\times V$.  Since  $L^2(\mu)$ is a Hilbert space and $\{f_n:n\in\N\}$ is bounded in $L^2(\mu)$ by \eqref{bound-f-n1}--\eqref{bound-f-n2},   there exists a subsequence $(n_k)_{k\in \N}$ and a $g\in L^2(\mu)$ such that $\{f_{n_k}\}$ converges to $g$ in $\sigma(L^2,L^2)$, see \cite[IV.4, Corollary~7]{Dunford}.  As in \eqref{eq_h}--\eqref{eq_f} it follows that for $m$-a.e.\ $v\in G$, $f(\cdot,v)$ is absolutely continuous  with derivative $g$. 

Step~\emph{(b3)}:
By  \eqref{inf-var},  $G\cup J$ is a $m$ null set,  and hence for $m$-a.e.\ $v$,  $f(\cdot,v)$  is absolutely continuous;  let $\dot f(\cdot,v)$ denote its derivative. Since $\dot f\in L^2(\mu)$, \eqref{gen_eq} follows and  we only need to show \eqref{fdot_int}.  Since for $m$-a.e.\ $v$, $f(\cdot,v)$ is  absolutely continuous with derivative $\dot f(\cdot,v)$ we have that $f_n\to \dot f$ $\lambda\otimes m$-a.e. By  continuity of $s\mapsto \xi_v(s)$, it follows that $\xi_v(f_n(s,v))\to \xi_v(\dot f(s,v))$ for $\lambda\otimes m$-a.e.\ $(s,v)$. Thus, by Fatou's Lemma and \eqref{xi_int},
\begin{align}\label{}
{}& \int_{V_k}\int_{\R} \xi_v(\dot f(s,v))\, ds\, m(dv) \leq \liminf_{n\to \infty} \int_{V_k} \int_{\R}\xi_v(f_n(s,v)) \, ds\, m(dv)<\infty,
\end{align}
which shows \eqref{fdot_int}. This completes the proof of \emph{(b)}. 

\emph{(a)}: In the general situation, define two (positive) L\'evy measures $\rho^1_v$ and $\rho^2_v$ by  
\begin{equation}\label{}
  \rho^1_v(dx)=\frac{1}{1\vee x^2}\,\rho_v(dx)\quad \text{and}\quad \rho^2_v=\rho_v-\rho^1_v, \qquad  v\in V,  
\end{equation}
and let   $X^1$ and $X^2$ be two independent processes defined as $X$ with $\rho$ replaced by $\rho^1=\{\rho^1_v:v\in V\}$ and  $\rho^2=\{\rho^2_v:v\in V\}$, respectively.  Since $X\dist X^1+X^2$, a symmetrization argument shows that $X^1$ is of finite variation. Moreover, since $\int_{\abs{x}>1} x^2\,\rho^1_v(dx)=\rho_v([-1,1]^c)<\infty$, Proposition~\ref{remark}\emph{(i)} shows that $\rho^1$ satisfies \eqref{eq:u0}, and hence  \eqref{gen_eq} and \eqref{eq:1suf} follow by \emph{(b)}. This completes the proof of \emph{(a)}. 

\emph{(c)}: Assume that $\rho$ satisfies \eqref{eq:u00}. This  yields the existence of a real constant $C_0>0$ such that for all  $u>0$ and $v\in V$ 
\begin{equation}\label{}
 \int_{\abs{ux}>1} \abs{xu}\,\rho_v(dx)\leq C_0 \int_{\R} (\abs{xu}^2\wedge 1)\,\rho_v(dx).
\end{equation}
Hence  for all $r>0$,
\begin{align}\label{eq:K_0:1}
 {}&  \int_{\{(x,v,s) \in \R\times V\times \R: \, \|xf(\cdot - s,v)\|_{BV[0,1]}>1\}}  \|xf(\cdot - s,v)\|_{BV[0,1]} \, \rho_v(dx)\, m(dv)  \, ds\\ \label{eq:K_0:2}{}& \qquad \leq C_0\int_{\R} \int_{V} \int_{\R} \left(1 \wedge \|xf(\cdot - s,v)\|_{BV[0,1]}^2 \right) \, \rho_v(dx)\, m(dv)\, ds  < \infty,
\end{align}
which by Lemma~\ref{int_ID} shows that $\E\norm{X^P}_{BV[0,1]}<\infty$. By arguing as above, \eqref{eq:Cf} follows. This completes the proof of \emph{(c)}. 
\end{proof}

\begin{remark}\label{remark-con}
In the proof of Theorem~\ref{necessary_fv} we  only used the assumption \eqref{inf-var} to conclude that $f(\cdot,v)$ is absolutely continuous for $m$-a.e.\ $v$. Thus if we know that $f(\cdot,v)$ is absolutely continuous for $m$-a.e.\  $v$ then Theorem~\ref{necessary_fv} remains valid even without the assumption \eqref{inf-var}.
 \end{remark}
 
 %%%%%%%%%%%%%%%%
\section{Proofs of Theorems~\ref{th:1} and \ref{th:2}}
\label{proof_01}
%%%%%%%%%%%%%%

\begin{proof}[Proof of Theorem~\ref{th:1}]
Recall that $\D$ denotes the set of dyadic numbers in $\R$. Consider $\R^{\D}$ as a locally convex separable linear metric space and consider $X_{\D}:=(X_t)_{t\in \D}$ as a random variable in $\R^{\D}$. 
For each  $N \in \N$, define
\begin{equation}\label{eq:HN}
H_N =  \Big\{ h \in \R^{\D}: \, \sup_{n\in \N}\sum_{i=1}^{2N2^{n}} \abs{h(r_{n,i}^N)- h(r_{n,i-1}^N)} < \infty \Big\},
\end{equation}
where $r_{n,i}^N= i2^{-n}-N$, and let $H= \bigcap_{N=1}^{\infty} H_N$. By \eqref{eq:fvX} 
\begin{equation}
\P\big(\|X\|_{BV[a,b]}< \infty \quad \text{for all } -\infty < a < b < \infty \big)=\P(X_{| \D}\in H).
\end{equation}
Let $\nu$ be the L\'evy measure of $X_{\D}$. We have
\begin{align}
\nu(H_N^c)  &= \int_{\R}\int_{\R\times V}  \1_{H_N^c}(xf(\cdot -s, v)) \, \rho_v(dx) \, m(dv) \, ds \\
& = \int_{\R}\int_{V} \rho_v(\R) \, \1_{H_N^c}(f(\cdot -s, v)) \, m(dv) \, ds \label{eq:Hc}
\end{align}
because $\rho_v(\{0\})=0$. By \eqref{eq:fvf} we also have
\begin{equation} \label{eq:stc-f}
\|f(\cdot-s, v)\|_{BV[-N,N]} = \sup_{n\in \N}\sum_{i=1}^{2N2^n} \Big{|}f(r_{n,i}^N -s, v) - f(r_{n,i-1}^N- s,v)\Big{|} \quad \lambda \otimes m\text{-a.e.}
\end{equation}
Consider the set 
\begin{equation}\label{}
A = \{v: \rho_v(\R)>0 \ \text{and} \ \|f(\cdot, v)\|_{BV[-M,M]}= \infty \ \text{for some } M \in \N\}.
\end{equation}
If $m(A)=0$ then $\nu(H_N^c)=0$ for every $N$, and so $\nu(H^c)=\lim_{N \to \infty}\nu(H_N^c)=0$. From \citet[Theorem 9]{Janssen}, we get $\P(X_{| \D}\in H)=0$ or 1. 

Suppose now that $m(A)>0$, so that $m(A_M)>0$ for some $M \in \N$, where
\begin{equation}\label{eq:0-1-c}
A_M = \{v: \rho_v(\R)>0 \ \text{and} \ \|f(\cdot, v)\|_{BV[-M,M]}= \infty \}.
\end{equation}
For every $N>M$ and  all $(s,v) \in [M-N, N-M] \times A_M$ 
we have 
\begin{equation} \label{}
\|f(\cdot-s,v)\|_{BV[-N,N]} \ge  \|f(\cdot, v)\|_{BV[-M,M]} = \infty,
\end{equation}
which combined with \eqref{eq:Hc} gives 
\begin{equation} \label{}
\nu(H_N^c) \ge 2(N-M) \int_{A_M} \rho_v(\R) \, m(dv).  
\end{equation}
Thus  $\nu(H^c)= \lim_{N \to \infty} \nu(H_N^c) = \infty$. By \citet[Theorem 10]{Janssen}, $\P(X_{| \D} \in H)=0$. This completes the proof.
\end{proof}

\begin{proof}[Proof of Theorem~\ref{th:2}]
Fix $a<b$ and define
\begin{equation}\label{}
  H= \Big\{\morf{h}{\D}{\R}:  \sup_{n\in \N} \ \sum_{i=1}^{k_n} 
|h(r_{n,i}) - h(r_{n, i-1})| < \infty\Big\},
\end{equation}
where $\{r_{n,i}\}$ are a dyadic partitions of $[a,b]$ such that $\max_{1\leq i\leq k_n} (r_{i,n}-r_{i-1,n})\to 0$ as $n\to \infty$.  As in \eqref{eq:Hc} we show that
\begin{align}\label{eq:Hc1}
\nu(H^c) &=   \int_{\R} \int_V  \rho_v(\R) \, \1_{H^c}(f(\cdot -s, v)) \,  m(dv) ds.
\end{align}
If \emph{(a)} holds then $\nu(H^c)=0$ and  the zero-one law holds by the same argument as in the previous theorem. 

Assume \emph{(b)}. Let $m(A_M)>0$ for some $M>0$, where $A_M$ is given by \eqref{eq:0-1-c}. Then
there exists a subinterval $[c,d] \subset [-M,M]$ with $d-c< (b-a)/2$ such that $m(B)>0$, where
 \begin{equation} \label{}
B:=\{v: \rho_v(\R) >0 \ \text{and} \ \|f(\cdot, v)\|_{BV[c,d]}= \infty\}.
\end{equation}
For all $(s,v) \in [a-c, b-d] \times B$ we have 
\begin{equation} \label{}
\|f(\cdot-s,v)\|_{BV[a,b]} \ge  \|f(\cdot, v)\|_{BV[-M,M]} = \infty,
\end{equation}
which combined with \eqref{eq:Hc1}  gives 
\begin{equation} \label{}
\nu(H^c) \ge \frac{b-a}{2} \int_{B} \rho_v(\R) \, m(dv) = \infty.  
\end{equation}
By the same argument as in the proof of Theorem \ref{th:1} we infer that the probability in \eqref{eq:0-1a} is zero. 
If $m(A_M)=0$ for all $M \in \N$, then $\nu(H^c)=0$. 
We conclude, as above, that the probability in \eqref{eq:0-1a} is 0 or 1.
\end{proof}
\medskip

Finally, let us note that  the methods of proofs of  Theorems~\ref{th:1} and \ref{th:2} will work if we replace the total variation norm $\norm{\,\cdot\,}_{BV[a,b]}$  by some wider class of  seminorms of sample paths. 

\bigskip\medskip\noindent
{\large \textbf{Acknowledgement} }
\medskip 

\noindent
We thank the referee for a  constructive and detailed report.

\bibliographystyle{chicago}

 %\bibliography{bibliografi}

\end{document}